\documentclass[10pt]{article}
\usepackage{psfrag,amsmath}
\usepackage{amssymb}
\usepackage{mathrsfs}
\usepackage{amsfonts}
\usepackage{amsmath}
\usepackage{mathrsfs,color}
\usepackage{epstopdf}
\usepackage{graphicx}
\usepackage{mathrsfs,amscd,amssymb,amsthm,amsmath,bm,url}

\makeatletter

\renewcommand{\@seccntformat}[1]{{\csname the#1\endcsname}{\normalsize .}\hspace{.5em}}
\makeatother

\def \[{\begin{equation}}
\def \]{\end{equation}}

\def \ex{{\rm ex}}

\def \rank{{\rm rank}}
\def \Tr{{\rm Tr}}
\newtheorem{thm}{Theorem}[section]

\newtheorem{claim}{Claim}
\newtheorem{defi}{Definition}
\newtheorem{fact}{Fact}
\newtheorem{lem}[thm]{Lemma}
\newtheorem{cor}[thm]{Corollary}

\newtheorem{pb}{Problem}

\newtheorem{conj}[thm]{Conjecture}
\newtheorem{remark}{Remark}

\newenvironment{wst}
{\setlength{\leftmargini}{1.5\parindent}
 \begin{itemize}
 \setlength{\itemsep}{-1.1mm}}
{\end{itemize}}
\begin{document}

\setlength{\baselineskip}{15pt}
\begin{center}{\Large \bf Adjacency eigenvalues of graphs without short odd cycles\footnote{Financially supported  by the National Natural Science Foundation of China (Grant Nos. 12171190, 11671164).
\\
\hspace*{5mm}{\it Email addresses}: lscmath@mail.ccnu.edu.cn (S.C. Li), \ wtsun2018@sina.com (W.T. Sun), \ ytyumath@sina.com
 (Y.T. Yu).}}
\vspace{4mm}

{\large Shuchao Li,\ \ Wanting Sun,\ \ Yuantian Yu}\vspace{2mm}

Faculty of Mathematics and Statistics,  Central China Normal University, Wuhan 430079, P.R. China
\end{center}

\noindent {\bf Abstract}:\ It is well known that spectral Tur\'{a}n type problem is one of the most classical {problems} in graph theory. In this paper, we consider the spectral Tur\'{a}n type problem. Let $G$ be a graph and let $\mathcal{G}$ be a set of graphs, we say $G$ is \textit{$\mathcal{G}$-free} if $G$ does not contain any element of $\mathcal{G}$ as a subgraph. Denote by $\lambda_1$ and $\lambda_2$ the largest and the second largest eigenvalues of the adjacency matrix $A(G)$ of $G,$ respectively. In this paper we focus on the characterization of graphs without short odd cycles according to the adjacency eigenvalues of the graphs. Firstly, an upper bound on $\lambda_1^{2k}+\lambda_2^{2k}$ of $n$-vertex $\{C_3,C_5,\ldots,C_{2k+1}\}$-free graphs is established, where $k$ is a positive integer. All the corresponding extremal graphs are identified. Secondly, a sufficient condition for non-bipartite graphs containing an odd cycle of length at most $2k+1$ in terms of its spectral radius is given. At last, we characterize the unique graph having the maximum spectral radius among the set of $n$-vertex non-bipartite graphs with odd girth at least $2k+3,$ which solves an open problem proposed by Lin, Ning and Wu [Eigenvalues and triangles in graphs, Combin. Probab. Comput. 30 (2) (2021) 258-270]. 

\vspace{2mm} \noindent{\it Keywords:}
Eigenvalue; Spectral radius; Odd cycle; Spectral Tur\'{a}n problem
\vspace{2mm}


\section{\normalsize Introduction}
We will begin with introducing the background information which will derive our main results. Our main results will also be given in this section.
\subsection{\normalsize Background and definitions}
In this paper, we consider only simple, undirected and finite graphs. Let $G=(V_{G},E_{G})$ be a graph, where $V_G$ is its vertex set and $E_G$ is its edge set. The \textit{order} of $G$ is the number $n=|V_G|$ of its vertices and its \textit{size} is the number $|E_G|$ of its edges. Denote by  $P_n,\,C_n,\,K_n$ and $K_{a,n-a}$ the path, the cycle, the complete graph and the complete bipartite graph on $n$ vertices, respectively. We say that two vertices $u$ and $v$ are \textit{adjacent} (or \textit{neighbors}) if they are joined by an edge. The {set of neighbors} of a vertex $u$ is denoted by $N_G(u)$ (or, $N(u)$ for short). The {\textit{degree}} $d_G(u)$ (or, $d(u)$ for short) of a vertex $u$ (in $G$) is the cardinality of $N_G(u)$. Then $\delta(G):=\min_{u\in V_G}d_G(u)$ is the \textit{minimum degree} of $G$. Unless otherwise stated, we follow the traditional notation and terminology; see \cite{0001,0009}.


Given a graph $G$, its \textit{adjacency matrix} $A(G)$ is an {$n\times n$\, $0$-$1$ matrix} whose $(i,j)$-entry is $1$ if and only if $ij\in E_G$. It is obvious that $A(G)$ is a real symmetric, nonnegative and irreducible matrix if $G$ is connected. {Hence,} its eigenvalues are real and can be given in non-increasing order as $\lambda_1(G)\geqslant \lambda_2(G)\geqslant \cdots \geqslant \lambda_n(G)$. In the whole context, when there is no scope for ambiguity, we always write $\lambda_i$ instead of $\lambda_i(G)$ for $1\leqslant i\leqslant n$.

{Denote by $s^+(G)$ (resp. $s^-(G)$) the sum of squares of positive (resp. negative) eigenvalues of $A(G).$} The largest modulus of an eigenvalue of $A(G)$ is called {the} \textit{spectral radius} of $G$. By the famous Perron-Frobenius theorem, we know that $\lambda_1$ is the spectral radius of $G$ and there exists a positive eigenvector ${\bf x} = (x_1, x_2, \ldots,x_n)^T$ of $A(G)$ corresponding to $\lambda_1$. It will be convenient to associate a labeling of vertices of $G$ (with respect to ${\bf x}$) in which $x_r$ is a label of the vertex $r$.

Let $G$  be a graph and let $\mathcal{G}$ be a set of graphs, we say that $G$ is \textit{$\mathcal{G}$-free} if it does not contain any graph in $\mathcal{G}$ as a subgraph. In particular, {if $\mathcal{G}=\{H\}$, then we also say that $G$ is $H$-free.} {The \textit{Tur\'{a}n number}, $\ex(n,\mathcal{G}),$ is} the maximum number of edges in a $\mathcal{G}$-free graph of order $n.$ To determine the exact value of $\ex(n,\mathcal{G})$ is a central problem of extremal graph theory, which is known as the Tur\'{a}n problem. However, there are only a few cases when the Tur\'{a}n number is known. For more details on this topic, the readers may be
referred to see \cite{19} for $\ex(n, K_{r+1}),$ \cite{02,03,7,21} for $\ex(n, C_{2l+1}),$ and \cite{4,6,8,9,17} for $\ex(n, C_{2l}).$

{Extremal  problems  involving  cycles have  been  considered  since the beginning of graph theory. It is interesting and challenging to study the structure or parameters of a graph if this graph has (no) short odd (or even) cycles. Ne\v{s}et\v{r}il and R\"{o}dl \cite{Ne} showed that graphs having no short odd cycles have the edge-partition property. Shearer \cite{She} studied the independence number of dense graphs with large odd girth. Gy\"{o}ri, Kostochka and {\L}uczak \cite{GKL} proved that graphs without short odd cycles are nearly bipartite. Very recently, Lin and Zeng \cite{LZ} studied the extremal problems on the bisection of graphs without short even cycles. We refer the reader to \cite{Dud,Har,FLZ,Low,Ped} for further results in this direction.

The question which interests us here is the question of what spectral condition guarantees a graph without short odd cycles.  We provide a characterization of all graphs which have no short odd cycles according to the eigenvalues of the adjacency matrices.}

It is well known that the spectra of graphs are a useful tool for us to investigate the graph parameters and graph structural properties. Particularly,  many researchers focus on bounding the spectral radius of a graph according to some of its classical parameters. In 1985, Brualdi and Hoffman \cite{Brua} showed that $\lambda_1(G)\leqslant k-1$ if $|E_G|\leqslant \frac{k(k-1)}{2}$ for some positive integer $k.$ {Wu and Elphick \cite{Wu} generalized this result as follows:
$$
  \sqrt{s^+(G)}\leqslant \frac{1}{2}\left(\sqrt{8|E_G|+1}-1\right).
$$
It also strengthens Stanley's inequality \cite{Stan} that:
$
  \lambda_1(G)\leqslant \frac{1}{2}(\sqrt{8|E_G|+1}-1).
$}
Another well-known upper bound for the spectral radius of a {connected graph}, due to Hong \cite{hong}, is that
$$
  \lambda_1(G)\leqslant \sqrt{2|E_G|-|V_G|+1},
$$
the equality holds only for complete graphs and star graphs. This bound has been strengthened by several authors. For example, Nikiforov \cite{Nik} proved that:
$$
  \lambda_1(G)\leqslant \frac{\delta(G)-1}{2}+\sqrt{2|E_G|-|V_G|\delta(G)+\frac{(1+\delta(G))^2}{4}},
$$
which is exact for various families of graphs, including regular graphs. {Furthermore, Elphick, Farber, Goldberg and Wocjan \cite{Elp} conjectured that if $G$ is a connected graph with order $n,$ then
$$
  \min\{s^+(G), s^-(G)\}\geqslant n-1,
$$
which has been proved for various classes of graphs, including bipartite, regular, complete $q$-partite, hyper-energetic and barbell graphs.} For more results, we refer the reader to \cite{1,2,3} and the references cited therein.

In 2010, Nikiforov \cite{Nik2} proposed a spectral version of extremal graph theory problem, which is also known as Brualdi-Solheid-Tur\'{a}n type problem, i.e., what is the maximal spectral radius of an $H$-free graph of order $n ?$ Over the past decade, much attention has been paid to the Brualdi-Solheid-Tur\'{a}n type problem. For more details, one may consult the references, such as for $H\cong K_r$ \cite{Nik3,Wilf}, $H\cong K_{s,t}$ \cite{Bab,Nik4,Nik3}, $H\cong P_k$ \cite{Nik2}, $H\cong C_4$ \cite{Nik3,ZW} and $H\cong C_6$ \cite{zhai}.

A \textit{$t$-walk} in a graph $G$ is an alternative sequence of vertices
and edges $v_1,e_1,v_2,e_2, \ldots,e_{t-1},v_t$ such that each edge $e_i$ is incident with $v_i$ and $v_{i+1}$ ($1\leqslant i\leqslant t-1$). The walk is \textit{closed} if $v_1$ coincides with $v_t$. Denote by $\alpha_t(G)$ and $\beta_t(G)$ the numbers of $t$-walks and closed $t$-walks in $G,$ respectively. Let $\omega=\omega(G)$ denote the clique number of $G.$ Wilf \cite{Wilf} determined that
$$
  \lambda_1(G)\leqslant \frac{\omega-1}{\omega}|V_G|=\frac{\omega-1}{\omega}\alpha_1(G).
$$
In \cite{Nik5}, Nikiforov generalized this result as follows: for each integer $s\geqslant 2,$
$$
  \lambda_1^s(G)\leqslant \frac{\omega-1}{\omega}\alpha_s(G).
$$
Notice that for the case $s=2,$ the above inequality implies the concise form of {Tur\'{a}n's} theorem. Therefore, Nikiforov's inequality sometimes is called the {\textit{spectral Tur\'{a}n theorem}}.

In 2007, Bollob\'{a}s and Nikiforov \cite{Boll} posed the following interesting conjecture, which gives an upper bound on $\lambda_1^2+\lambda_2^2$ for $K_{r+1}$-free graphs.
\begin{conj}[\cite{Boll}]\label{conj1}
Let $G$ be a $K_{r+1}$-free graph of order at least $r+1$ with size $m.$ Then
$$
  \lambda_1^2+\lambda_2^2\leqslant \frac{r-1}{r}2m.
$$
\end{conj}
Recently, Lin, Ning and Wu \cite{lin} confirmed this conjecture in the case $r=2,$ and they also characterized all the extremal graphs. Motivated by their work, our first result {(i.e., Theorem \ref{thm1.1})} establishes an upper bound on $\lambda_1^{2k}+\lambda_2^{2k}$ among the set of $\{C_3,C_5,\ldots,C_{2k+1}\}$-free graphs for each positive integer $k,$ which generalizes the result of Lin, Ning and Wu \cite[Theorem 1.2]{lin}.

On the other hand, we also notice that Mantel's theorem is a quintessential result in extremal graph theory, which studies the maximum number of edges over all $C_3$-free graphs:
\begin{thm}[Mantel's Theorem]
   Every graph of order $n$ and size greater than $\lfloor\frac{n^2}{4}\rfloor$ contains a triangle.
\end{thm}
Erd\H{o}s \cite[Exercise 12.2.7]{Bon} focused on this problem among the set of non-bipartite $C_3$-free graphs. Nosal \cite{2} established the spectral version of Mantel's theorem: each $C_3$-free graph $G$ satisfies $\lambda_1(G)\leqslant \sqrt{|E_G|}.$ Lin, Ning and Wu \cite{lin} proved two spectral analogues results, i.e., each non-bipartite graph $G$ of order $n$ and size $m$ contains $C_3$ as a subgraph if one of the following holds: (i) $\lambda_1\geqslant \sqrt{m-1}$ and $G\not\cong C_5\cup (n-5)K_1;$ (ii) $\lambda_1\geqslant \lambda_1(S(K_{\lfloor(n-1)/2\rfloor,\lceil(n-1)/2\rceil}))$ and $G\not\cong S(K_{\lfloor(n-1)/2\rfloor,\lceil(n-1)/2\rceil}),$ where $S(K_{\lfloor(n-1)/2\rfloor,\lceil(n-1)/2\rceil})$ is obtained from $K_{\lfloor(n-1)/2\rfloor,\lceil(n-1)/2\rceil}$ by subdividing an edge. Along this line, our second main result {(i.e., Theorem \ref{thm1.3})} gives a sufficient condition for non-bipartite graphs containing an odd cycle of length at most $2k+1.$

Now, let us recall the extension of the famous Erd\H{o}s's theorem:
\begin{thm}\label{thm1.04}
Let $k$ be a positive integer. Assume that $G$ is a $\{C_3,C_5,\ldots,C_{2k+1}\}$-free graph with order $n.$ If $G$ is non-bipartite, then
\[\label{eq10}
      |E_G|\leqslant \left(\frac{n-2k+1}{2}\right)^2+2k-1.
\]
\end{thm}
{Lin, Ning and Wu \cite{lin} provided a non-spectral proof of Theorem \ref{thm1.04}. {Motivated by this nice result,} they proposed the following problem:}
\begin{pb}\label{pb1}
Let $k$ be a positive integer. How can we characterize all the graphs among the set of $\{C_3,C_5,\ldots,C_{2k+1}\}$-free non-bipartite graphs with order $n$ achieving the maximum spectral radius?
\end{pb}
{In view of the result in \cite[Theorem 1.4]{lin}, this problem was settled for the case $k=1$.} Our last main result {(i.e., Theorem \ref{thm1.4})} completely solves this open problem for each positive integer $k$.

\subsection{\normalsize Main results}
In this subsection, we give our main results. {Recall that $P_n$ is a path of order $n$. Let $t\geqslant 2$ be a positive integer. Then $tP_n$ denotes the disjoint union of $t$ copies of $P_n.$} For two vertex-disjoint graphs $G_1$ and $G_2$, the \textit{join} $G_1\vee G_2$ is the graph obtained by joining every vertex of $G_1$ with every vertex of $G_2$. For a given graph $H,$ a  \textit{blow-up} of $H$ is a graph obtained from $H$ by replacing each vertex $v$ of $H$ with a stable set $I_v,$ in which we add all edges between $I_u$ and $I_v$ if $uv\in E_H.$

Our first result establishes an upper bound on $\lambda_1^{2k}+\lambda_2^{2k}$ among the set of $\{C_3,C_5,\ldots,C_{2k+1}\}$-free graphs for each positive integer $k,$ and all the corresponding extremal graphs are characterized.
\begin{thm}\label{thm1.1}
Let $k$ be a positive integer and let $G$ be a graph of order $n\, (\geqslant 2k+1).$ If $G$ is $\{C_3,C_5,\ldots,C_{2k+1}\}$-free, then
\[\label{eq:1}
  \lambda_1^{2k}+\lambda_2^{2k}\leqslant {\frac{\Tr(A^{2k}(G))}{2}}=\frac{\beta_{2k}(G)}{2},
\]
{where $\Tr(\cdot)$ denotes the trace of a matrix.} The equality holds if and only if $G$ is a blow-up of $H$ with
$$
  H\in \{P_2\cup K_1,2P_2\cup K_1,P_4\cup K_1,P_5\cup K_1\}.
$$
\end{thm}

The subsequent result gives a sufficient condition for non-bipartite graphs containing an odd cycle of length at most $2k+1.$
\begin{thm}\label{thm1.3}
Let $k$ be a positive integer and let $G$ be a non-bipartite graph. 
{If $G$ is $\{C_3,C_5,\ldots,C_{2k+1}\}$-free, then
$$
  \lambda_1^{2k}\leqslant {\frac{\Tr(A^{2k}(G))}{2}-\left(2\cos\frac{\pi}{k+2}\right)^{2k}}.
$$
The equality holds if and only if $k=1$ and $G$ is isomorphic to $C_{5}$ together with some isolated vertices.} 
\end{thm}

Our last main result gives a complete solution for Problem \ref{pb1}. Notice that for each $\{C_3,C_5,\ldots,C_{2k+1}\}$-free non-bipartite graph $G,$ one has $|V_G|\geqslant 2k+3$ and $G\cong C_{2k+3}$ if $|V_G|= 2k+3.$ So, it suffices to consider the case $n\geqslant 2k+4$ in Problem~\ref{pb1}. Let $R_k(K_{s,t})$ be a graph obtained from a complete bipartite graph $K_{s,t}$ by replacing one of its edges with {$P_{2k+1}.$}

\begin{thm}\label{thm1.4}
Let $k$ be a positive integer and let $G$ be a $\{C_3,C_5,\ldots,C_{2k+1}\}$-free non-bipartite graph with order $n\geqslant 2k+4$. Then
$$
\lambda_1(G)\leqslant \lambda_1(R_k(K_{\lfloor(n-2k+1)/2\rfloor,\lceil(n-2k+1)/2\rceil})).
$$
The equality holds if and only if $G\cong R_k(K_{\lfloor(n-2k+1)/2\rfloor,\lceil(n-2k+1)/2\rceil}).$
\end{thm}

The remainder of this paper is organized as follows: In Section~2, we give some preliminary results, which will be used in the subsequent sections. In section~3, we give the proofs of Theorems \ref{thm1.1} and \ref{thm1.3}. In Section~4, we present the proof of Theorem \ref{thm1.4}. In the last section, we give some further discussions.
\section{\normalsize Preliminaries}\setcounter{equation}{0}
In this section, we describe some preliminary results, which play an important role in the subsequent sections.
Let ${\bf r}=(r_1,r_2,\ldots,r_n)^T$ and ${\bf s}=(s_1,s_2,\ldots,s_n)^T$ be two column vectors in $\mathbb{R}^n.$ Now, we rearrange the elements of ${\bf r}$ and ${\bf s}$ in non-increasing orders as follows: $r_{[1]}\geqslant r_{[2]}\geqslant \cdots \geqslant r_{[n]}$ and $s_{[1]}\geqslant s_{[2]}\geqslant \cdots \geqslant s_{[n]}.$ Recall the definition of ``a vector is weakly majorized by the other one" as follows:
\begin{defi}
Let ${\bf r}=(r_1,r_2,\ldots,r_n)^T$ and ${\bf s}=(s_1,s_2,\ldots,s_n)^T$ be in $\mathbb{R}^n.$ We say that ${\bf s}$ is weakly majorized by ${\bf r}$ if
$$
  \sum_{i=1}^k s_{[i]}\leqslant \sum_{i=1}^k r_{[i]}
$$
for $k=1,2,\ldots,n$ and denote it by ${\bf s}\prec_w{\bf r}.$
\end{defi}
The following lemma {will be used to prove} Theorem \ref{thm1.1}. 
\begin{lem}[\cite{lin}]\label{lem2.1}
Let ${\bf r}=(r_1,r_2,\ldots,r_n)^T$ and ${\bf s}=(s_1,s_2,\ldots,s_n)^T$ be in $\mathbb{R}^n$ such that
$$
  r_1\geqslant r_2\geqslant \cdots \geqslant r_n\geqslant 0\ \text{and}\ s_1\geqslant s_2\geqslant \cdots \geqslant s_n\geqslant 0.
$$
If $p>1$ is a real number and ${\bf s}\prec_w{\bf r},$ then $\parallel{\bf s}\parallel_p\leqslant \parallel{\bf r}\parallel_p$ and the equality holds if and only if ${\bf r}={\bf s},$ where $\parallel{\bf r}\parallel_p=(\sum_{i=1}^n|r_i|^p)^{\frac{1}{p}}.$
\end{lem}
The next lemma is one of the fundamental results in spectral graph theory.
\begin{lem}[\cite{0001}]\label{lem3.4}
Let $G$ be a connected graph and let $H$ be a proper subgraph of $G.$ Then $\lambda_1(H)<\lambda_1(G).$ In addition, if $H$ is an induced subgraph of $G,$ then $\lambda_2(H)\leqslant \lambda_2(G).$
\end{lem}
Let $G$ be a graph. The \textit{rank} of $G,$ written as $\rank(G),$ is defined to be the rank of its adjacency matrix $A(G).$ The subsequent lemma characterizes the graphs with small rank.
\begin{lem}[\cite{Ob}]\label{lem2.2}
Let $G$ be a graph on $n$ vertices.
\begin{wst}
\item[{\rm (i)}] $\rank(G)=2$ if and only if $G$ is a blow-up of $P_2\cup K_1.$
\item[{\rm (ii)}] If $G$ is a bipartite graph, then $\rank(G)=4$ if and only if $G$ is a blow-up of $H,$ where $H\in\{2P_2\cup K_1,P_4\cup K_1,P_5\cup K_1\}.$
\end{wst}
\end{lem}
Let $T_{a,b,c}$ denote a $T$-\textit{shaped tree} defined as a tree with a single vertex $u$ of degree $3$ such that $T_{a,b,c}-u=P_a\cup P_b\cup P_c\, (1\leqslant a\leqslant b\leqslant c).$ 
\begin{lem}[\cite{Guo}]\label{lem3.3}
The eigenvalues of $A(T_{1,1,n-3})$ are $0$ and $2\cos\frac{(2k-1)\pi}{2n-2},\, k=1,\ldots,n-1.$
\end{lem}
Let $H$ be a real matrix, whose columns and rows are indexed by $U=\{1,2,\ldots,n\}.$ Assume that $\pi:=U_1\cup U_2\cup \ldots\cup U_t$ is a partition of $U$. Then $H$ can be partitioned based on $\pi$, i.e.,
\begin{equation*}
    H=\left(
        \begin{array}{ccc}
          H_{11} & \cdots & H_{1t} \\
          \vdots & \ddots & \vdots \\
          H_{t1} & \cdots & H_{tt} \\
        \end{array}
      \right),
\end{equation*}
where $H_{ij}$ denotes the submatrix of $H$, indexed by the rows and columns of $U_i$ and $U_j$ respectively. Let $\pi_{ij}$ be the average row sum of $H_{ij}$ for $1\leqslant i,j \leqslant t.$ Usually, the matrix $H^{\pi}=(\pi_{ij})$ is called the \textit{quotient matrix} of $H$. Moreover, if the row sum of $H_{ij}$ is constant for $1\leqslant i,j \leqslant t$, then we call $\pi$ an \textit{equitable partition}.
\begin{lem}[\cite{0009}]\label{lem3.1}
Let $H$ be a real matrix with an equitable partition $\pi$, and let $H^{\pi}$ be the corresponding quotient matrix. Then every eigenvalue of $H^{\pi}$ is an eigenvalue of $H$. In addition, if $H=A(G)$ for some graph $G,$ then the spectral radius of $G$ is equal to the {largest eigenvalue} of $H^{\pi}$.
\end{lem}
{\section{\normalsize Proofs of Theorems \ref{thm1.1} and \ref{thm1.3}}\setcounter{equation}{0}
In this section, we give the proofs of Theorems \ref{thm1.1} and \ref{thm1.3}. First, we prove Theorem \ref{thm1.1}, which establishes an upper bound on $\lambda_1^{2k}+\lambda_2^{2k}$ for an $n$-vertex $\{C_3,C_5,\ldots,C_{2k+1}\}$-free graph, and all the corresponding extremal graphs are characterized.} 

\begin{proof}[\bf Proof of Theorem \ref{thm1.1}]\
Let $p^+$ and $p^-$ denote the numbers (including the multiplicities) of positive and negative eigenvalues of the adjacency matrix $A(G)$, respectively. Put
$$
s_k^+:=\lambda_1^{2k}+\lambda_2^{2k}+\cdots+\lambda_{p^+}^{2k}\ \text{and}\ s_k^-:=\lambda_{n-p^-+1}^{2k}+\lambda_{n-p^-+2}^{2k}+\cdots+\lambda_{n}^{2k}.
$$
It is well known that $s_k^+ +s_k^-{=\Tr(A^{2k}(G))}=\beta_{2k}(G)$. Notice that $G$ is $C_3$-free. {Hence,} $G\not\cong K_n$ and so $\lambda_2(G)\geqslant 0$  (see \cite[Lemma 5]{10}).

Suppose to the contrary that $\lambda_1^{2k}+\lambda_2^{2k}> {\frac{\Tr(A^{2k}(G))}{2}}.$ Therefore,
$$
  \lambda_1^{2k}+\lambda_2^{2k}>\frac{s_k^++s_k^-}{2}.
$$
It follows that
\[\label{eq:2.1}
  \lambda_1^{2k}+\lambda_2^{2k}\geqslant 2(\lambda_1^{2k}+\lambda_2^{2k})-s_k^+>s_k^-\geqslant 0.
\]
Let
$$
  {\bf x}=(\lambda_1^{2k},\lambda_2^{2k},0,\ldots,0)^T\ \text{and}\ {\bf y}=(\lambda_{n-p^-+1}^{2k},\lambda_{n-p^-+2}^{2k},\ldots,\lambda_{n}^{2k})^T
$$
be two nonnegative vectors in $\mathbb{R}^{p^-}.$ Together with \eqref{eq:2.1}, we obtain ${\bf y}\neq {\bf x}$ and ${\bf y}\prec_w{\bf x}.$ Applying Lemma~\ref{lem2.1} with $p=\frac{2k+1}{2k}$ yields that
$$
  \parallel {\bf y}\parallel_{\frac{2k+1}{2k}}< \parallel {\bf x}\parallel_{\frac{2k+1}{2k}},\ \ \text{i.e.},\ \ (\parallel {\bf y}\parallel_{\frac{2k+1}{2k}})^{\frac{2k+1}{2k}}< (\parallel {\bf x}\parallel_{\frac{2k+1}{2k}})^{\frac{2k+1}{2k}},
$$
which is equivalent to
$$
  |\lambda_{n-p^-+1}^{2k+1}|+|\lambda_{n-p^-+2}^{2k+1}|+\cdots+|\lambda_{n}^{2k+1}|<\lambda_1^{2k+1}+\lambda_2^{2k+1}.
$$
Therefore,
\begin{align*}
  \beta_{2k+1}(G)&={\Tr(A^{2k+1}(G))}=\lambda_1^{2k+1}+\lambda_2^{2k+1}+\cdots+\lambda_{p^+}^{2k+1}+\lambda_{n-p^-+1}^{2k+1}+\lambda_{n-p^-+2}^{2k+1}+\cdots+\lambda_{n}^{2k+1}\\[5pt]
  &\geqslant \lambda_1^{2k+1}+\lambda_2^{2k+1}+\lambda_{n-p^-+1}^{2k+1}+\lambda_{n-p^-+2}^{2k+1}+\cdots+\lambda_{n}^{2k+1}>0,
\end{align*}
which implies that there exists a closed $(2k+1)$-walk in $G.$ It is straightforward to check that $G$ contains at least one odd cycle $C_{2i+1}$ with $1\leqslant i\leqslant k,$ a contradiction. {Hence,} the inequality in \eqref{eq:1} holds.

Now, we assume that the equality in \eqref{eq:1} holds, i.e., $\lambda_1^{2k}+\lambda_2^{2k}={\frac{\Tr(A^{2k}(G))}{2}}=\frac{s_k^++s_k^-}{2}.$ Clearly, $\lambda_1^{2k}+\lambda_2^{2k}\geqslant s_k^-\geqslant 0.$ Therefore, ${\bf y}\prec_w{\bf x}.$ Applying Lemma~\ref{lem2.1} with $p=\frac{2k+1}{2k}$, one has
$$
  \parallel {\bf y}\parallel_{\frac{2k+1}{2k}}\leqslant \parallel {\bf x}\parallel_{\frac{2k+1}{2k}},\ \ \text{i.e.},\ \ (\parallel {\bf y}\parallel_{\frac{2k+1}{2k}})^{\frac{2k+1}{2k}}\leqslant (\parallel {\bf x}\parallel_{\frac{2k+1}{2k}})^{\frac{2k+1}{2k}}.
$$
It follows that
$$
  |\lambda_{n-p^-+1}^{2k+1}|+|\lambda_{n-p^-+2}^{2k+1}|+\cdots+|\lambda_{n}^{2k+1}|\leqslant \lambda_1^{2k+1}+\lambda_2^{2k+1}.
$$
Notice that $G$ is an $n$-vertex $\{C_3,C_5,\ldots,C_{2k+1}\}$-free graph. Hence,
\begin{align*}
  0=\beta_{2k+1}(G)&={\Tr(A^{2k+1}(G))}=\lambda_1^{2k+1}+\lambda_2^{2k+1}+\cdots+\lambda_{p^+}^{2k+1}+\lambda_{n-p^-+1}^{2k+1}+\lambda_{n-p^-+2}^{2k+1}+\cdots+\lambda_{n}^{2k+1}\\[5pt]
  &\geqslant \lambda_1^{2k+1}+\lambda_2^{2k+1}+\lambda_{n-p^-+1}^{2k+1}+\lambda_{n-p^-+2}^{2k+1}+\cdots+\lambda_{n}^{2k+1}\geqslant 0.
\end{align*}
Thus,
$$
  \lambda_1^{2k+1}+\lambda_2^{2k+1}=-(\lambda_{n-p^-+1}^{2k+1}+\lambda_{n-p^-+2}^{2k+1}+\cdots+\lambda_{n}^{2k+1}),
$$
which is equivalent to
$$
  (\parallel {\bf x}\parallel_{\frac{2k+1}{2k}})^{\frac{2k+1}{2k}}=(\parallel {\bf y}\parallel_{\frac{2k+1}{2k}})^{\frac{2k+1}{2k}},\ \ \text{i.e.},\ \  \parallel {\bf x}\parallel_{\frac{2k+1}{2k}}=\parallel {\bf y}\parallel_{\frac{2k+1}{2k}}.
$$
Applying Lemma \ref{lem2.1} again yields that ${\bf x}={\bf y},$ i.e.,
$$
  \lambda_1^{2k}=\lambda_n^{2k},\ \lambda_2^{2k}=\lambda_{n-1}^{2k}\ \text{and}\ \lambda^2_{n-p^-+1}=\cdots=\lambda^2_{n-3}=0.
$$
That is to say,
$$
\lambda_1=-\lambda_n,\ \lambda_2=-\lambda_{n-1}\ \text{and}\ \lambda_{n-p^-+1}=\cdots=\lambda_{n-3}=0.
$$
Notice that $\Tr(A(G))=\sum_{i=1}^n\lambda_i=0.$ {Hence,} $\lambda_3=\cdots=\lambda_{n-3}=0$ and therefore $G$ is a bipartite graph.

If $\lambda_2=0,$ then $\rank(G)=2.$ Together with Lemma \ref{lem2.2}(i), one obtains that $G$ is a blow-up of $P_2\cup K_1.$ If $\lambda_2\neq 0,$ then $\rank(G)=4.$ Based on Lemma \ref{lem2.2}(ii), we know that $G$ is a blow-up of $H$ with $H\in \{2P_2\cup K_1,P_4\cup K_1,P_5\cup K_1\},$ as desired.

Conversely, if $G$ is a blow-up of $H$ with $H\in \{P_2\cup K_1,2P_2\cup K_1,P_4\cup K_1,P_5\cup K_1\},$ then $G$ is bipartite. Together with Lemma~\ref{lem2.2}, we see that
$$
  \lambda_1(G)=-\lambda_n(G),\ \lambda_2(G)=-\lambda_{n-1}(G)\ \text{and}\ \lambda_{3}(G)=\cdots=\lambda_{n-3}(G)=0.
$$
Therefore, $\lambda_1^{2k}+\lambda_2^{2k}={\frac{\Tr(A^{2k}(G))}{2}}.$ This completes the proof.
\end{proof}
The following result is an immediate consequence of Theorem \ref{thm1.1}.
\begin{cor}\label{thm1.2}
Let $G$ be a graph. If $\lambda_1^{2k}(G)\geqslant {\frac{\Tr(A^{2k}(G))}{2}},$ then $G$ contains an odd cycle with length at most $2k+1$ unless $G$ is a blow-up of $P_2\cup K_1.$
\end{cor}
\begin{proof}
Suppose to the contrary that $G$ is $\{C_3,C_5,\ldots,C_{2k+1}\}$-free and $G$ is not isomorphic to the blow-up of $P_2\cup K_1.$ Since  $\lambda_1^{2k}\geqslant {\frac{\Tr(A^{2k}(G))}{2}},$ one has $\lambda_1^{2k}+\lambda_2^{2k}\geqslant {\frac{\Tr(A^{2k}(G))}{2}}.$ Together with Theorem~\ref{thm1.1}, we obtain that $\lambda_1^{2k}+\lambda_2^{2k}= {\frac{\Tr(A^{2k}(G))}{2}}$ and therefore $G$ is a blow-up of $H$ with $H\in\{2P_2\cup K_1,P_4\cup K_1,P_5\cup K_1\}.$ In addition, we know that $\lambda_2=0$ and $G$ is bipartite. On the other hand, in view of Lemma \ref{lem2.2}(ii), we have $\rank(G)=4,$ which implies that $\lambda_2>0,$ a contradiction.
\end{proof}
\begin{remark}{\rm
In Corollary \ref{thm1.2}, put $k=1$, then it is just the main result \cite[Theorem 2(i)]{Nik6}. Hence, Corollary \ref{thm1.2} gives a new and simple proof for \cite[Theorem 2(i)]{Nik6}.}
\end{remark}
\vspace{2mm}
{Now, we are ready to prove Theorem \ref{thm1.3}.} Recall that $T_{a,b,c}$ denotes a $T$-shaped tree. It is routine to check that
\[\label{eq:2}
\lambda_1(T_{1,2,2})=2\cos\frac{\pi}{12},\, \lambda_1(T_{1,2,3})=2\cos\frac{\pi}{18}\ \text{and}\   2\cos\frac{\pi}{31}>\lambda_1(T_{1,2,4})>2\cos\frac{\pi}{30}.
\]

\begin{proof}[\bf Proof of Theorem \ref{thm1.3}]\
{Notice that Lin, Ning and Wu \cite[Theorem 1.3]{lin} confirmed this result for $k=1.$ So, we only need to consider the case that $k\geqslant 2.$}  {In order to prove this result, it suffices to show that if $k\geqslant 2$ and $G$ is a non-bipartite graph with
\[\label{eq:001}
  \lambda_1^{2k}\geqslant {\frac{\Tr(A^{2k}(G))}{2}-\left(2\cos\frac{\pi}{k+2}\right)^{2k}},
\]
then $G$ contains at least one graph among $\{C_3,C_5,\ldots,C_{2k+1}\}.$} 

Suppose that $G$ is a $\{C_3,C_5,\ldots,C_{2k+1}\}$-free non-bipartite graph. Hence, $|V_G|\geqslant 2k+3.$ At first we show the following claim.
\begin{claim}\label{claim1}
$\lambda_2(G)< 2\cos\frac{\pi}{k+2}.$
\end{claim}
\begin{proof}[\bf Proof of Claim \ref{claim1}]
If this is not true, then $\lambda_2(G)\geqslant 2\cos\frac{\pi}{k+2}.$ Hence, {\eqref{eq:001} implies}
$$
  \lambda_1^{2k}+\lambda_2^{2k}\geqslant {\frac{\Tr(A^{2k}(G))}{2}}-\left(2\cos\frac{\pi}{k+2}\right)^{2k}+\left(2\cos\frac{\pi}{k+2}\right)^{2k}={\frac{\Tr(A^{2k}(G))}{2}}.
$$
Together with Theorem \ref{thm1.1}, one has $\lambda_1^{2k}+\lambda_2^{2k}= {\frac{\Tr(A^{2k}(G))}{2}}$ and so $G$ is a blow-up of $H$ with $H\in\{P_2\cup K_1,2P_2\cup K_1,P_4\cup K_1,P_5\cup K_1\}.$ Obviously, $G$ is bipartite, a contradiction. Thus, $\lambda_2(G)<2\cos\frac{\pi}{k+2}.$
\end{proof}

We proceed by considering whether $G$ is connected or not.
We firstly consider the case that $G$ is connected. Let $C_s=u_1u_2\ldots u_su_1$ be a shortest odd cycle of $G.$ It follows that $C_s$ is an induced cycle of $G$ and $s\geqslant 2k+3.$ It is well known to us that $\lambda_2(C_s)=2\cos \frac{2\pi}{s}.$ If $s\geqslant 2k+5,$ then by Lemma \ref{lem3.4} one has
$$
  \lambda_2(G)\geqslant \lambda_2(C_s)=2\cos \frac{2\pi}{s}>2\cos \frac{2\pi}{2k+4}=2\cos\frac{\pi}{k+2},
$$
a contradiction to Claim \ref{claim1}. Therefore, $s=2k+3.$ {If $|V_G|=2k+3,$ then $G\cong C_{2k+3}.$ Notice that each closed $2k$-walk starting with the vertex $u_i$ in $C_{2k+3}$ is determined by the set of times in the walk in which we move ``forward" (i.e., from $u_j$ to $u_{j+1}$). Since there are $k$ ``forward" and $k$ ``backward" moves, one has $\Tr(A^{2k}(C_{2k+3}))=\beta_{2k}(C_{2k+3})=(2k+3){{2k}\choose {k}}.$ Next, in order to obtain a contradiction to \eqref{eq:001}, we are to prove that if $G\cong C_{2k+3}$ with $k\geqslant 2,$ then
\[\label{eq:003}
  2^{2k}<\frac{(2k+3){{2k}\choose {k}}}{2}-\left(2\cos\frac{\pi}{k+2}\right)^{2k}.
\]
If $2\leqslant k\leqslant 9,$ then by Mathematica 9.0 we get \eqref{eq:003} immediately. For $k\geqslant 10,$ it suffices to show
$$
  2^{2k}<\frac{(2k+3){{2k}\choose {k}}}{2}-2^{2k},
$$
which is equivalent to
\[\label{eq:002}
  4<\frac{(2k+3){{2k}\choose {k}}}{2^{2k}}.
\]

Now, we use induction on $k$ to prove \eqref{eq:002}. If $k=10,$ then by a direct calculation we obtain \eqref{eq:002} immediately. Assume that \eqref{eq:002} is true for each integer no more than $k\,(\geqslant 10).$ By induction, it is routine to check that
$$
  \frac{(2(k+1)+3){{2(k+1)}\choose {k+1}}}{2^{2(k+1)}}=\frac{(2k+5)(2k+2)(2k+1)}{4(2k+3)(k+1)^2}\frac{(2k+3){{2k}\choose {k}}}{2^{2k}}>\frac{2(2k+5)(2k+1)}{(2k+3)(k+1)}>4.
$$
Hence \eqref{eq:002} holds for all integers $k\geqslant 10,$ as desired.} So, in what follows, we assume that $|V_G|\geqslant 2k+4.$

Let $v$ be a vertex in $N(u_i)\setminus V_{C_s}$ for some $i\,(1\leqslant i\leqslant 2k+3).$ Recall that $G$ is $\{C_3,C_5,\ldots,C_{2k+1}\}$-free. Then $N(v)\cap V_{C_s}\subseteq \{u_i,u_{i+2}\},$ taking subscripts modulo $2k+3.$ {For each positive integer $k,$ denote by $H_1^k$ the graph obtained from the cycle $C_{2k+3}$ and the vertex $v$ by adding an edge $vu_i$, and let $H_2^k=H_1^k+vu_{i+2}.$ It follows that either  $H_1^k$ or $H_2^k$ is an induced subgraph of $G.$ 

Note that $H_1^k-u_{i+2}\cong T_{1,1,2k}$ and $H_2^k-u_{i+2}\cong T_{1,1,2k}.$ By Lemmas \ref{lem3.4} and \ref{lem3.3}, one has for $k\geqslant 2,$
$$
  \lambda_2(G)\geqslant\min\{\lambda_2(H_1^k),\lambda_2(H_2^k)\}\geqslant \lambda_2(T_{1,1,2k})=2\cos\frac{3\pi}{4k+4}\geqslant 2\cos\frac{\pi}{k+2},
$$
which contradicts Claim \ref{claim1}.}

Now, we consider the case that $G$ is disconnected. {Choose  a connected component, say $G_1,$ of $G$ such that $\lambda_1(G)=\lambda_1(G_1).$ Notice that $\Tr(A^{2k}(G))\geqslant\Tr(A^{2k}(G_1)).$ In view of \eqref{eq:001}, we get
$$
\lambda_1^{2k}(G_1)=\lambda_1^{2k}(G)\geqslant {\frac{\Tr(A^{2k}(G))}{2}-\left(2\cos\frac{\pi}{k+2}\right)^{2k}}\geqslant \frac{\Tr(A^{2k}(G_1))}{2}-\left(2\cos\frac{\pi}{k+2}\right)^{2k}.
$$
Applying the discussion in the connected case to $G_1$ yields that $G_1$ contains an odd cycle of length at most $2k+1.$ It follows that $G$ contains at least one graph among $\{C_3,C_5,\ldots,C_{2k+1}\},$ which contradicts the choice of $G.$}

This completes the proof.
\end{proof}

Let $G$ be a graph of size $m$ with maximum degree $\Delta,$ Chen and Qian \cite{CQ} showed that the number of closed walks $\beta_l(G)$ with $l\geqslant 3$ in $G$ satisfies $\beta_l(G)\leqslant 2m\Delta^{l-2}$ with equality if and only if $l$ is even and each component of $G$ is a complete bipartite graph $K_{\Delta,\Delta}.$ Together with Theorem \ref{thm1.3}, we obtain the following corollary immediately.
\begin{cor}\label{cor001}
Let $k\geqslant 2$ be a positive integer and let $G$ be a non-bipartite graph with size $m$ and maximum degree $\Delta.$ If $G$ is $\{C_3,C_5,\ldots,C_{2k+1}\}$-free, then
$$
  \lambda_1<\sqrt[2k]{m\Delta^{2k-2}-\left(2\cos\frac{\pi}{k+2}\right)^{2k}}.
$$
\end{cor}

\begin{remark}{\rm
In Theorem \ref{thm1.3}, if $k=1$, then $\lambda_1\leqslant \sqrt{m-1}$ and the equality holds if and only if $G$ is isomorphic to $C_5$ together with some isolated vertices. Note that it is just the main result \cite[Theorem 1.3]{lin}.

Now assume that $k\geqslant 2$ is a positive integer and $G$ is a non-bipartite $\{C_3,C_5,\ldots,C_{2k+1}\}$-free graph of size $m$ with maximum degree $\Delta,$ where $\Delta\leqslant \sqrt[2k-2]{\frac{(m-1)^k}{m}}.$ Then by a direct calculation one has
\begin{eqnarray*}
  \lambda_1(G)&\leqslant& \sqrt[2k]{\frac{\Tr(A^{2k}(G))}{2}-\left(2\cos\frac{\pi}{k+2}\right)^{2k}}\\
              &\leqslant& \sqrt[2k]{m\Delta^{2k-2}-\left(2\cos\frac{\pi}{k+2}\right)^{2k}}\\
              &<&\sqrt[2k]{m\Delta^{2k-2}}\\
              &\leqslant& \sqrt{m-1}.
\end{eqnarray*}
Therefore, when $k\geqslant2$ the upper bounds in both Theorem \ref{thm1.3} and Corollary~\ref{cor001} are better than that of \cite[Theorem 1.3]{lin}}.
\end{remark}
{\section{\normalsize Proof of Theorem \ref{thm1.4}}\setcounter{equation}{0}
In this section, we give the proof of Theorem \ref{thm1.4}. Before doing this, we need the following lemma.} Recall that $R_k(K_{s,t})$ denotes the graph obtained by replacing an edge of the complete bipartite graph $K_{s,t}$ with {$P_{2k+1}.$}
\begin{lem}\label{lem3.2}
Let $s,t,k$ be three positive integers with $t\geqslant s\geqslant 2$. If $t-s\geqslant 2,$ then
$$
  \lambda_1(R_k(K_{s+1,t-1}))> \lambda_1(R_k(K_{s,t})).
$$
\end{lem}
\begin{proof}
It is easy to see that $C_{2k+3}$ is a proper subgraph of both $R_k(K_{s+1,t-1})$ and $R_k(K_{s,t}).$ By Lemma~\ref{lem3.4}, one has
$$
  \lambda_1(R_k(K_{s+1,t-1}))>\lambda_1(C_{2k+3})=2\ \text{and}\ \lambda_1(R_k(K_{s,t}))>\lambda_1(C_{2k+3})=2.
$$
Let $X\cup Y$ be the bipartition of {$V_{K_{s,t}}$} with $|X|=s$ and $|Y|=t.$ Without loss of generality, we assume that the graph $R_k(K_{s,t})$ is obtained from $K_{s,t}$ by replacing the edge $uv$ with the path $uu_1u_2\ldots u_{2k-1}v,$ where $u\in X$ and $v\in Y.$ It is straightforward to check that $\pi:=(X\setminus \{u\})\cup (Y\setminus \{v\})\cup \{u\}\cup \{u_1\}\cup\{u_2\}\cup\ldots\cup \{u_{2k-1}\}\cup \{v\}$ is an equitable partition of $R_k(K_{s,t})$ {with respect to $A(R_k(K_{s,t})).$} Hence, the quotient matrix of $A(R_k(K_{s,t}))$ corresponding to the partition $\pi$ can be {given by}
\begin{equation*}
    (A(R_k(K_{s,t})))^{\pi}=\left[
       \begin{array}{cccccccc}
         0 & t-1 & 0 & 0 & 0 & \cdots & 0 & 1 \\
         s-1 & 0 & 1 & 0 & 0 & \cdots & 0 & 0 \\
         0 & t-1 & 0 & 1 & 0 & \cdots & 0 & 0 \\
         0 & 0 & 1 & 0 & 1 & \cdots & 0 & 0 \\
         0 & 0 & 0 & 1 & 0 & \cdots & 0 & 0 \\
         \vdots & \vdots & \vdots & \vdots & \vdots & \ddots & \vdots & \vdots\\
         0 & 0 & 0 & 0 & 0 & \cdots & 0 & 1 \\
         s-1 & 0 & 0 & 0 & 0 & \cdots & 1 & 0 \\
       \end{array}
     \right],
\end{equation*}
which is a $(2k+3)$-by-$(2k+3)$ matrix.

Put $f_k(\lambda,s,t):=\det(\lambda I-(A(R_k(K_{s,t})))^{\pi})$ and $g_n(\lambda):=\det(\lambda I-A(P_{n}))$. Then evaluating $\det(\lambda I-(A(R_k(K_{s,t})))^{\pi})$ by Laplace expansion in terms of the first two columns gives us
\begin{align*}
  f_k(\lambda,s,t)=&\left|
                     \begin{array}{cc}
                       \lambda & -(t-1) \\
                       -(s-1) & \lambda \\
                     \end{array}
                   \right|g_{2k+1}(\lambda)+\left|
                     \begin{array}{cc}
                       \lambda & -(t-1) \\
                       0 & -(t-1) \\
                     \end{array}
                   \right|g_{2k}(\lambda)+\left|
                     \begin{array}{cc}
                       \lambda & -(t-1) \\
                       -(s-1) & 0 \\
                     \end{array}
                   \right|\\[5pt]
                   &-\left|
                     \begin{array}{cc}
                     -(s-1) & \lambda \\
                     0 & -(t-1) \\
                     \end{array}
                   \right|-\left|
                     \begin{array}{cc}
                     -(s-1) & \lambda \\
                     -(s-1) &  0\\
                     \end{array}
                   \right|g_{2k}(\lambda)-\left|
                     \begin{array}{cc}
                       0 & -(t-1) \\
                       -(s-1) & 0\\
                     \end{array}
                   \right|g_{2k-1}(\lambda)\\[5pt]
                    =&(\lambda^2-(s-1)(t-1))g_{2k+1}(\lambda)-(s+t-2)\lambda g_{2k}(\lambda)-2(s-1)(t-1)+(s-1)(t-1)g_{2k-1}(\lambda).
\end{align*}
Therefore,
$$
  f_k(\lambda,s+1,t-1)=(\lambda^2-s(t-2))g_{2k+1}(\lambda)-(s+t-2)\lambda g_{2k}(\lambda)-2s(t-2)+s(t-2)g_{2k-1}(\lambda).
$$
{{Notice that $g_{n+1}(\lambda)=\lambda g_{n}(\lambda)-g_{n-1}(\lambda).$} Hence,}
$$
  f_k(\lambda,s+1,t-1)-f_k(\lambda,s,t)=(s-t+1)(g_{2k+1}(\lambda)-g_{2k-1}(\lambda)+2)=(s-t+1)(\lambda g_{2k}(\lambda)-2g_{2k-1}(\lambda)+2).
$$

Let $\lambda_0$ be an arbitrary real number greater than $2.$ In what follows, we use induction on $n$ to prove $g_{n+1}(\lambda_0)-g_{n}(\lambda_0)>0$ for any integer $n\geqslant 1.$ Clearly the statement holds for $n=1.$ Assume that the statement is true for $n-1.$ 
By a direct calculation, one has
$$
  g_{n+1}(\lambda_0)-g_{n}(\lambda_0)=\lambda_0 g_{n}(\lambda_0)-g_{n-1}(\lambda_0)-g_{n}(\lambda_0)=(\lambda_0-1)g_{n}(\lambda_0)-g_{n-1}(\lambda_0).
$$
By the inductive assumption, we have
$$
  g_{n+1}(\lambda_0)-g_{n}(\lambda_0)> (\lambda_0-2)g_{n-1}(\lambda_0)>0,
$$
the last inequality holds since the largest eigenvalue of {$A(P_{n-1})$ is $2\cos{\frac{\pi}{n}}<2.$} Thus, $g_{n+1}(\lambda_0)-g_{n}(\lambda_0)>0$ holds for each positive integer $n.$
It follows that {$g_{n+1}(\lambda)-g_{n}(\lambda)>0$} for each integer $n\geqslant 1$ and each real number $\lambda>2.$ Therefore, $\lambda g_{2k}(\lambda)-2g_{2k-1}(\lambda)+2>0$  for any $\lambda>2.$

Recall that $t-s\geqslant 2.$ {Thus,}  if $\lambda>2,$ then
$$
  f_k(\lambda,s+1,t-1)-f_k(\lambda,s,t)<0.
$$
That is to say, the largest eigenvalue of $(A(R_k(K_{s+1,t-1})))^{\pi}$ is greater than that of $(A(R_k(K_{s,t})))^{\pi}.$ Together with Lemma~\ref{lem3.1}, one has $\lambda_1(R_k(K_{s+1,t-1}))> \lambda_1(R_k(K_{s,t})).$ 
This completes the proof.
\end{proof}
Next, we close this section by giving the proof of Theorem \ref{thm1.4}.
\begin{proof}[\bf Proof of Theorem \ref{thm1.4}]\
Assume that $G$ is a $\{C_3,C_5,\ldots,C_{2k+1}\}$-free non-bipartite graph of order $n$ having the maximum spectral radius. In order to complete the proof, it suffices to show that $G$ is isomorphic to $R_k(K_{\lfloor(n-2k+1)/2\rfloor,\lceil(n-2k+1)/2\rceil}).$

Firstly, we show that $G$ is connected. Otherwise, $G$ contains at least two connected components, say $G'$ and $G''.$ Let $\tilde{G}$ be a graph obtained from $G$ by connecting a vertex in $G'$ and some other vertex in $G''$ by an edge. Clearly, $\tilde{G}$ is a $\{C_3,C_5,\ldots,C_{2k+1}\}$-free non-bipartite graph of order $n$. Furthermore, by Lemma \ref{lem3.4} one has $\lambda_1(\tilde{G})>\lambda_1(G),$ which contradicts the choice of $G$.

Let ${\bf x}=(x_1,\ldots,x_n)^T$ be the unit positive eigenvector of $A(G)$ corresponding to $\lambda_1(G)$ and let
$$
  x_\triangleleft=\max\{x_i:1\leqslant i\leqslant n\},
$$
which corresponds to the vertex $u$ in $G$.

Notice that $G$ is non-bipartite and $\{C_3,C_5,\ldots,C_{2k+1}\}$-free. Without loss of generality, we may assume that $C=u_1u_2\ldots u_su_1$ is a shortest odd cycle of $G$ with $s\geqslant 2k+3.$ Clearly, $C$ is an induced subgraph of $G.$ We shall {characterize the structure of $G$} by considering the following facts.
\begin{fact}\label{c2}
$s=2k+3.$
\end{fact}
\begin{proof}[\bf Proof of Fact \ref{c2}]
Suppose to the contrary that $s\geqslant 2k+5.$ Recall that $C$ is an induced subgraph of $G.$ It follows that $u_1u_4\not\in E_G.$ It is easy to see that $G+u_1u_4$ is non-bipartite and $\lambda_1(G+u_1u_4)>\lambda_1(G)$ (based on Lemma~\ref{lem3.4}).  By the choice of $G,$ we know that $G+u_1u_4$ contains {a shortest} odd cycle $C_{2i+1}^1$ of {length $2i+1\,(1\leqslant i\leqslant k)$ with} $u_1u_4\in E_{C_{2i+1}^1}.$ Assume that $C_{2i+1}^1=u_1v_2\ldots v_{2i}u_4u_1.$ {Hence,} there exist at least two distinct $(u_1,u_4)$-paths in $G,$ one of which is $u_1u_2u_3u_4$ and the other is $u_1v_2\ldots v_{2i}u_4.$

Put $B:=\{u_2,u_{3}\}\cap V_{C_{2i+1}^1}.$ Clearly, $0\leqslant |B|\leqslant 2.$ If $|B|=0,$ then $u_1u_2u_3u_4 v_{2i}\ldots v_2 u_1$ forms an odd cycle of length $2i+3$ in $G,$ a contradiction. {If $|B|=1,$ then assume that $u_2\in B$ and $v_j=u_2.$ It is routine to check that one of $u_1v_2\ldots v_{j-1}u_2u_1$ and $u_2v_{j+1}\ldots v_{2i}u_4u_3u_2$ is an odd cycle of length no more than $2i+1$ in $G,$ a contradiction. If $|B|=2,$ then assume that $u_2=v_j$ and $u_3=v_l.$ If $j<l,$ then $j$ and $l$ are even. Otherwise, without loss of generality, suppose that $j$ is odd. Then $u_1v_2\ldots v_{j-1}u_2u_1$ is an odd cycle of length $j$ in $G,$ a contradiction. Therefore,  $u_2v_{j+1}\ldots v_{l-1}u_3u_2$ is an odd cycle of length $l-j+1$ in $G,$ a contradiction. If $l<j,$ by a similar discussion we obtain that $j$ and $l$ are odd. Therefore, $u_3v_{l+1}\ldots v_{j-1}u_2u_3$ is an odd cycle of length $j-l+1$ in $G,$ a contradiction.} Therefore, $s=2k+3.$ This completes the proof of Fact \ref{c2}.
\end{proof}
\begin{fact}\label{c3}
For each vertex $w\in V_G\setminus (N(u)\cup V_C),$ one has $N(w)=N(u).$
\end{fact}
\begin{proof}[\bf Proof of Fact \ref{c3}]\
If $V_G\setminus (N(u)\cup V_C)=\emptyset,$ then we are done. Now, we assume that $V_G\setminus (N(u)\cup V_C)\neq \emptyset.$ Suppose to the contrary that there exists a vertex $w\in V_G\setminus (N(u)\cup V_C)$ such that $N(w)\neq N(u).$ Let
$$
G_1:=G-\{wz:z\in N(w)\}+\{wz:z\in N(u)\}.
$$
Obviously, $G_1$ is non-bipartite.

Now, we show that $G_1$ is $\{C_3,C_5,\ldots,C_{2k+1}\}$-free. Otherwise, $G_1$ contains {a shortest odd cycle $C_{2i+1}^2$ of length $2i+1\,(1\leqslant i\leqslant k)$ with $w\in V_{C_{2i+1}^2}.$ It follows that $C_{2i+1}^2$ is an induced subgraph of $G_1.$ By the structure of $G_1,$ we obtain  $u\not\in V_{C_{2i+1}^2}.$} Then there exists an odd cycle with vertices $(V_{C_{2i+1}^2}\setminus \{w\})\cup \{u\}$ in $G,$ a contradiction. 

On {the} other hand,
\[\label{eq:3.1}
  \lambda_1(G_1)-\lambda_1(G)\geqslant {\bf x}^T(A(G_1)-A(G)){\bf x}= 2x_w\left(\sum_{z\in N_G(u)}x_z-\sum_{z\in N_G(w)}x_z\right)\geqslant 0.
\]
If $N_G(w)\subsetneq N_G(u),$ then $\lambda_1(G_1)>\lambda_1(G),$ a contradiction. It follows that
$$
  N_G(w)\setminus N_G(u)\neq \emptyset\ \ \text{and}\ \ N_G(u)\setminus N_G(w)\neq \emptyset.
$$
Recall that $\lambda_1(G)\geqslant \lambda_1(G_1).$ {Hence,} all inequalities in \eqref{eq:3.1} must be equalities. Therefore, ${\bf x}$ is an eigenvector of $A(G_1)$ corresponding $\lambda_1(G_1).$ Let $v$ be in $N_G(u)\setminus N_G(w).$ Then
$$
  \lambda_1(G)x_v=\sum_{z\in N_G(v)}x_z<\sum_{z\in N_{G_1}(v)}x_z=\lambda_1(G_1)x_v,
$$
i.e., $\lambda_1(G)<\lambda_1(G_1),$ which contradicts the choice of $G$. This completes the proof of Fact \ref{c3}.
\end{proof}
\begin{fact}\label{c4}
Let $x_{u'}=\max\{x_v:v\in N(u)\setminus V_C\}$ correspond to the vertex $u'$ in $N(u)\setminus V_C$. Then for any vertex $w$ other than $u'$ in $N(u)\setminus V_C,$ one has
$$
  N(w)\cap V_C=N({u'})\cap V_C.
$$
\end{fact}
\begin{proof}[\bf Proof of Fact \ref{c4}]
Notice that $G$ is $C_3$-free. Then $N_G(u)$ is a stable set of $G.$ For convenience, let
$$
  X:= N_G(u)\setminus V_C, \ \ \ \ Y:=\{y: y\in V_G\setminus V_C,\, N_G(u)=N_G(y)\}.
$$
In view of Fact \ref{c3}, we obtain that for all $v\in X,$
$$
  \lambda_1(G)x_v=\sum_{z\in Y}x_z+\sum_{z\in N_G(v)\cap V_C}x_z.
$$
It follows that for all $v\in X,$
\[\label{eq:4.2}
  \sum_{z\in N_G({u'})\cap V_C}x_z\geqslant \sum_{z\in N_G(v)\cap V_C}x_z.
\]

Suppose to the contrary that there exists a vertex $w\in X$ such that $N(w)\cap V_C\neq N({u'})\cap V_C.$ Then let
$$
  G_2:=G-\{wz:z\in N_G(w)\cap V_C\}+\{wz:z\in N_G({u'})\cap V_C\}.
$$
Clearly, $G_2$ is non-bipartite. Now, we show that $G_2$ is $\{C_3,C_5,\ldots, C_{2k+1}\}$-free. Otherwise, $G_2$ contains {a shortest} odd cycle $C_{2i+1}^3$ of {length $2i+1\,(i\in \{1,2,\ldots, k\})$ with} $ww'\in E_{C_{2i+1}^3}$ for some $w'\in N_G({u'})\cap V_C.$ Without loss of generality, we assume that $N_{G_2}(w)\cap V_{C_{2i+1}^3}=\{w',w''\}.$

If $w''\not\in V_C,$ then $ww''\in E_G.$ {Hence,} $w''\not\in N_G(u).$ Otherwise, $uww''u$ is a cycle of length $3$ in $G,$  a contradiction. In view of Fact \ref{c3}, one has $N_G(w'')=N_G(u)$ and therefore $w''\in N_G({u'}).$ If $w''\in V_C,$ based on the structure of $G_2,$ we have $w''\in N_G({u'}).$ That is to say, in both cases we have $w''\in N_G({u'}).$  {Recall that $C_{2i+1}^3$ is an induced subgraph of $G_2.$ Therefore, $u'\not\in V_{C_{2i+1}^3}.$} Then there exists an odd cycle with vertices $(V_{C_{2i+1}^3}\setminus \{w\})\cup \{u'\}$ in $G,$ a contradiction. It follows that $G_2$ is $\{C_3,C_5,\ldots, C_{2k+1}\}$-free.

On the other hand, it is straightforward to check that
\begin{align*}
  \lambda_1(G_2)-\lambda_1(G)\geqslant &\ {\bf x}^T(A(G_2)-A(G)){\bf x}\notag\\
                             = &\ 2x_w\left(\sum_{z\in N_G(u')\cap V_C}x_z-\sum_{z\in N_G({w})\cap V_C}x_z\right)\notag\\
                             \geqslant &\ 0,
\end{align*}
the last inequality follows by \eqref{eq:4.2}. By a similar discussion as the proof of Fact \ref{c3}, we obtain that $\lambda_1(G_2)>\lambda_1(G),$ which contradicts the choice of $G.$

This completes the proof of Fact \ref{c4}.
\end{proof}
Recall that
$$
  X=N_G(u)\setminus V_C\ \ \text{and}\ \ Y=\{y: y\in V_G\setminus V_C,\, N_G(u)=N_G(y)\}.
$$
Then in view of Fact \ref{c3}, one has $V_G=X\cup Y\cup V_C$ and $X\cap Y=\emptyset.$ For convenience, we put $|X|=a$ and $|Y|=b.$

In order to characterize the structure of $G$, we proceed by distinguishing the following four possible cases.

{\bf Case 1.} \ $u\not\in V_C$ and $X\neq \emptyset.$ In this case, note that $u\in Y$ and therefore $Y\neq \emptyset.$ Together with Facts \ref{c3} and \ref{c4}, one obtains $V_G\setminus V_C=X\cup Y$ and $G-C\cong K_{a,b}.$ Furthermore, applying Facts \ref{c3} and \ref{c4} again yields that each pair of distinct vertices in $X$ (resp. $Y$) share the same neighbors in $C.$ Denote by $N_C(S)$ the set of vertices in $C$ in which each is adjacent to some vertices of $S$ in $G,$ where $S\subseteq V_G\setminus V_C.$ Let $d_C(S)=|N_C(S)|.$ We proceed by showing the following claim.
\begin{claim}\label{c50}
$N_C(X)=\{u_i,u_{i+2}\}$ and $N_C(Y)=\{u_{i+1},u_{i+3}\}$ for some $i\in \{1,2,\ldots,2k+3\},$ where the subscripts are computed modulo $2k+3.$
\end{claim}
\begin{proof}[\bf Proof of Claim \ref{c50}]
Recall that $G$ is $\{C_3,C_5,\ldots,C_{2k+1}\}$-free and the shortest odd cycle of $G$ has length $2k+3$ (based on Fact~\ref{c2}). Let $v$ be an arbitrary vertex in $V_G\setminus V_C.$ Then $v$ has at most $2$ neighbors in $C$ and $N_C(v)\subseteq \{u_j,u_{j+2}\}$ for some $j\in\{1,\ldots,2k+3\}.$ It follows that $d_C(X)\leqslant 2$ and $d_C(Y)\leqslant 2.$

We are to prove that $d_C(X)=d_C(Y)=2.$ Since $G$ is connected, one has $\max\{d_C(X),d_C(Y)\}\geqslant 1.$ We may assume, without loss of generality, that $d_C(X)\geqslant 1.$ Suppose that $d_C(X)=1.$ By symmetry, we assume that $N_C(X)=\{u_1\}.$ Since $G$ is $\{C_3,C_5,\ldots,C_{2k+1}\}$-free, one has $N_C(Y)\subseteq \{u_2,u_4,u_{2k+1},u_{2k+3}\}.$ {Without loss of generality, we assume that $|N_C(Y)\cap \{u_2,u_4\}|\geqslant |N_C(Y)\cap \{u_{2k+1},u_{2k+3}\}|.$}

Note that $u_3,\,u_{2k+2}\not\in N_C(Y)$ if $k\geqslant 2$, and at least one of $\{u_3,u_4\}$ does not in $N_C(Y)$ if $k=1.$ Then we may assume that $u_3\not\in N_C(Y).$ Now, we construct a new graph $G_3$ by adding all edges between $\{u_3\}$ and $X.$ It is easy to see that $G_3$ is non-bipartite and $\lambda_1(G_3)>\lambda_1(G)$ (based on Lemma \ref{lem3.4}). By the choice of $G,$ we obtain that $G_3$ contains {a shortest} odd cycle $C_{2i+1}^4$ of {length $2i+1\,(1\leqslant i\leqslant k)$ with} $xu_3\in E_{C_{2i+1}^4}$ for some $x\in X.$ Note that $x$ and $u_3$ have no common neighbor in $G_3.$ Hence, $k\geqslant i\geqslant 2.$ Thus $u_3,u_{2k+2}\not\in N_C(Y).$

Next, we show that $xu_1\in E_{C_{2i+1}^4}.$ Suppose to the contrary that $xu_1\not\in E_{C_{2i+1}^4}.$ Then there exists a vertex $y$ in $Y$ such that $xy\in E_{C_{2i+1}^4}$ and so $G_3$ contains a $(y,u_3)$-path, say $yy_2\ldots y_{2i-1}u_3.$ {Clearly, $y_{2i-1}\not\in X.$ Otherwise, $yy_{2i-1}\in E_G$ and so $yy_2\ldots y_{2i-1}y$ is an odd cycle with length $2i-1$ of $G,$ a contradiction. It follows that $yy_2\ldots y_{2i-1}u_3$ is also a $(y,u_3)$-path in $G.$} {If $u_2\in N_C(Y)$, then there are at least two distinct paths that connect $y$ and $u_3$ in $G,$ one of which is $P^1:=yu_2u_3$ and the other is $P^2:=yy_2\ldots y_{2i-1}u_3.$  If $u_2\not\in V_{P^2},$ then $yy_2\ldots y_{2i-1}u_3u_2y$ is a cycle of length $2i+1\,(\leqslant 2k+1)$ in $G,$ a contradiction. If $u_2\in V_{P^2},$ then assume that $y_j=u_2.$ {Hence,} one of $yy_2\ldots y_{j-1}u_2y$ and $u_2y_{j+1}\ldots y_{2i-1}u_3u_2$ is an odd cycle of length no more than $2i-1$ in $G,$ a contradiction. Thus, $u_2\not\in N_C(Y). $ By a similar discussion, we can get $u_4\not\in N_C(Y).$ Recall that $|N_C(Y)\cap \{u_2,u_4\}|\geqslant |N_C(Y)\cap \{u_{2k+1},u_{2k+3}\}|,$ which implies $N_C(Y)= \emptyset,$ a contradiction.} {Thus,} $xu_1\in E_{C_{2i+1}^4}.$ Therefore, $G$ contains an odd cycle with vertex set $(V_{C_{2i+1}^4}\setminus \{x\})\cup \{u_2\},$ a contradiction. So, we obtain that $d_C(X)=2.$


Without loss of generality, assume that $N_C(X)=\{u_1,u_3\}.$ By symmetry, we know that
$$
  N_C(Y)\in \{\emptyset,\{u_2\},\{u_4\},\{u_2,u_4\}\}.
$$
If $N_C(Y)\subsetneq \{u_2,u_4\},$ then we construct a graph $G_4$ by adding all edges between $Y$ and $\{u_2,u_4\}.$ Clearly, $G_4\cong R_k(K_{a+2,b+2})$ and therefore $G$ is non-bipartite and $\{C_3,C_5,\ldots,C_{2k+1}\}$-free. In addition, by Lemma \ref{lem3.4}, one has $\lambda_1(G_4)>\lambda_1(G),$ a contradiction.

Therefore, $N_C(X)=\{u_1,u_3\}$ and $N_C(Y)=\{u_2,u_4\}$ in $G.$ This completes the proof of Claim \ref{c50}.
\end{proof}
Together with Fact \ref{c2} and Claim \ref{c50}, we deduce that $G\cong R_k(K_{a+2,b+2}).$ In view of Lemma \ref{lem3.2} and the choice of $G$, we have
$$
  G\cong R_k(K_{\lfloor(n-2k+1)/2\rfloor,\lceil(n-2k+1)/2\rceil}),
$$
as desired.

{\bf Case 2.}\ $u\not\in V_C$ and $X=\emptyset.$ In this case, $u\in Y$ and therefore $Y\neq \emptyset.$ Since $G$ is connected, one has $u_i\in N_C(Y)$ for some $1\leqslant i\leqslant 2k+3.$ Note that $d_C(Y)\leqslant 2.$ Based on the choice of $G$, one has $d_C(Y)=2.$ Otherwise, we can construct a graph $G_5$ by adding all edges between $Y$ and $\{u_{i+2}\}.$ Then $G_5\cong R_k(K_{2,b+2}).$ Clearly, $G_5$ is non-bipartite and $\{C_3,C_5,\ldots,C_{2k+1}\}$-free. By Lemma~\ref{lem3.4}, one has $\lambda_1(G_5)>\lambda_1(G),$ a contradiction. It follows that $G\cong R_k(K_{2,n-2k-1}).$ Applying Lemma \ref{lem3.2} yields
$$
  \lambda_1(G)= \lambda_1(R_k(K_{2,n-2k-1}))\leqslant \lambda_1(R_k(K_{\lfloor(n-2k+1)/2\rfloor,\lceil(n-2k+1)/2\rceil})),
$$
with equality if and only if $n=2k+4.$ Therefore, $G\cong R_k(K_{2,3})$ and $n=2k+4.$

{\bf Case 3.}\ $u\in V_C$ and $X\neq \emptyset.$ Without loss of generality, assume that $u=u_1.$ Based on Facts \ref{c3} and \ref{c4}, one has $N_C(Y)=\{u_2,u_{2k+3}\}$ and $u_1\in N_C(X).$ According to the choice of $G$, by a similar discussion as that in Case 2, we obtain that $u_3\in N_C(X)$ or $u_{2k+2}\in N_C(X).$ 
{Hence,} $G\cong R_k(K_{a+2,b+2}).$ In view of Lemma~\ref{lem3.2} and the choice of $G$, we have
$$
  G\cong R_k(K_{\lfloor(n-2k+1)/2\rfloor,\lceil(n-2k+1)/2\rceil}),
$$
as desired.

{\bf Case 4.}\ $u\in V_C$ and $X=\emptyset.$ Without loss of generality, assume that $u=u_1.$ Based on Fact \ref{c3}, one has $N_C(Y)=\{u_2,u_{2k+3}\}.$ Then  $G\cong R_k(K_{2,n-2k-1}).$  Applying Lemma \ref{lem3.2} yields
$$
  \lambda_1(G)= \lambda_1(R_k(K_{2,n-2k-1}))\leqslant \lambda_1(R_k(K_{\lfloor(n-2k+1)/2\rfloor,\lceil(n-2k+1)/2\rceil})),
$$
with equality if and only if $n=2k+4.$ Therefore, $G\cong R_k(K_{2,3})$ and $n=2k+4.$

This completes the proof.
\end{proof}
{It is well known that $\lambda_1(G)\geqslant \frac{2|E_G|}{|V_G|}$ and the equality holds if and only if $G$ is regular. So, the following result is an immediate consequence of Theorem \ref{thm1.4}.
\begin{cor}\label{cor1.4}
Let $k$ be a positive integer and $G$ be a $\{C_3,C_5,\ldots,C_{2k+1}\}$-free non-bipartite graph with order $n\geqslant 2k+4$. Then
\[\label{eq:11}
|E_G|< \frac{n}{2}\lambda_1(R_k(K_{\lfloor(n-2k+1)/2\rfloor,\lceil(n-2k+1)/2\rceil})).
\]
\end{cor}
\begin{remark}{\rm
Based on Lemma \ref{lem3.1}, it is straightforward to check that if $n$ is odd and $k\in\{1,2,3\}$, then
$$
  \lambda_1(R_k(K_{\lfloor(n-2k+1)/2\rfloor,\lceil(n-2k+1)/2\rceil}))>\left(\frac{n-2k+1}{2}\right)^2+2k-1.
$$
Thus, the upper bound in \eqref{eq10} is better than that of \eqref{eq:11} if $n$ is odd and $k\in\{1,2,3\}.$ But we do not know which one is better in the general case, since the exact value of $\lambda_1(R_k(K_{\lfloor(n-2k+1)/2\rfloor,\lceil(n-2k+1)/2\rceil}))$ can not be computed directly.}
\end{remark}}
\section{\normalsize Further discussions}
In this paper, we focus on graphs without short odd cycles. Firstly, we establish an upper bound on $\lambda_1^{2k}(G)+\lambda_2^{2k}(G)$ if $G$ is a $\{C_3,C_5,\ldots,C_{2k+1}\}$-free graph, and all the corresponding extremal graphs are characterized. It is interesting to see that our result is a natural generalization of one main results of Lin, Ning and Wu \cite{lin}. {Therefore, we confirm the Bollob\'{a}s-Nikiforov's conjecture (i.e., Conjecture \ref{conj1}) for the case $r=2$ in a new way.}

We also give a sufficient condition for non-bipartite graphs containing at least one graph in $\{C_3,C_5,\ldots,C_{2k+1}\}.$ In addition, we completely solve an open problem of Lin, Ning and Wu \cite{lin}, which determines the unique graph among the set of $n$-vertex non-bipartite  graphs with odd girth at least $2k+3$ having the maximum spectral radius.

The spectral Tur\'{a}n type problem is very interesting. It attracts more and more researchers' attention. In fact, some challenging problems on this topic are worthwhile studying.

{Notice that Nikiforov \cite{Nik} characterized the graph among the set of $n$-vertex $C_3$-free graphs having the maximum spectral radius; he also characterized the graph among the set of $n$-vertex $C_4$-free graphs having the maximum spectral radius (see \cite{Nik6}). Very recently, Zhai, Lin and Shu \cite{ZLS} characterized the graph among the set of $n$-vertex $C_5$-free (or $C_6$-free) graphs having the maximum spectral radius. We summarize them as follows.

For convenience, we use the notation $S_{n,k}$ to denote the graph obtained by the join of $K_k$ and $(n-k)K_1$, i.e., $S_{n,k}=K_k\vee (n-k)K_1$, where $n$ and $k$ be two positive integers.
\begin{thm}\label{thm6}
Let $G$ be an $H$-free graph with $m$ edges. The following assertions hold.
\begin{wst}
  \item[{\rm (a)}] If $H\cong C_3,$ then $\lambda_1(G)\leqslant \sqrt{m}$ with equality if and only if $G$ is a complete bipartite graph {\rm(\cite{Nik});}
  \item[{\rm (b)}] If $H\cong C_4$ with $m\geqslant 10,$ then $\lambda_1(G)\leqslant \sqrt{m}$ with equality if and only if $G$ is a star {\rm(\cite{Nik6});}
  \item[{\rm (c)}] If $H\cong C_5$ with $m\geqslant 8,$ or $H\cong C_6$ with $m\geqslant 22,$ then $\lambda_1(G)\leqslant \frac{1+\sqrt{4m-3}}{2}$ with equality if and only if $G\cong S_{\frac{m+3}{2},2}$ {\rm(\cite{ZLS})}.
\end{wst}
\end{thm}

In view of Theorem \ref{thm6}(a)-(b), we know that if $\lambda_1(G)\geqslant \sqrt{m},$ then $G$ contains $C_3$ and $C_4$ unless $G$ is a complete bipartite graph. By Theorem \ref{thm6}(c), one obtains that if $\lambda_1(G)\geqslant \frac{1+\sqrt{4m-3}}{2},$ then $G$ contains $C_t$ for every $t\leqslant 6$ unless $G\cong S_{\frac{m+3}{2},2}.$

Inspired by Theorem \ref{thm6}, Zhai, Lin and Shu \cite{ZLS} proposed the following conjecture (see, \cite[Conjecture 5.1]{ZLS}), which gives a more general spectral characterization of graphs containing cycles with consecutive lengths.
\begin{conj}\label{conj1.6}
Let $G$ be a graph of sufficiently large size $m$ without isolated vertices and let $k \geqslant 2$ be an integer. If
$$
\lambda_1(G)\geqslant \frac{k-1+\sqrt{4m-k^2+1}}{2},
$$
then $G$ contains $C_l$ for every $l\leqslant 2k+2$ unless $G\cong S_{\frac{m}{k}+\frac{k+1}{2},k}.$
\end{conj}}

Elphick, Farber, Goldberg and Wocjan \cite{Elp} studied the sum of squares of positive eigenvalues of a graph. It is natural and interesting for us to study the sum of even powers of positive eigenvalues of a graph. Thus, Theorem \ref{thm1.1} can be ported to a more general form.
\begin{pb}
Let $k$ be a positive integer and $G$ be a graph of order $n\ (\geqslant 2k+1).$ If $G$ is $\{C_3,C_5,\ldots,C_{2k+1}\}$-free, then how to determine the upper bound of $s_k^+,$ where $s_k^+$ is defined in the proof of Theorem~$\ref{thm1.1}$.
\end{pb}
We will do above open problems in the near future.

\section*{\normalsize Declaration of competing interest}
The authors declare that they have no known competing financial interests or personal relationships that could have
appeared to influence the work reported in this paper

\section*{\normalsize Acknowledgement}
We take this opportunity to thank the anonymous referees for their careful reading of the manuscript and suggestions which have immensely helped us in getting the article to its present form.

\end{document}